\documentclass[11pt]{article}

\usepackage[utf8]{inputenc}
\usepackage[T1]{fontenc}
\usepackage{lmodern}
\usepackage{amsmath,amssymb,amsfonts}
\usepackage{array}
\usepackage{booktabs}
\usepackage{longtable}
\usepackage{geometry}
\usepackage[hidelinks]{hyperref}
\usepackage{microtype}

\newcommand{\Z}{\mathbb Z}
\newcommand{\C}{\mathbb C}
\newcommand{\F}{\mathbb F}
\newcommand{\Ncal}{\mathcal N}
\newcommand{\Lcal}{\mathcal L}
\newcommand{\degtwo}{\deg_2}

\title{On the existence of Newman and Littlewood multiples for certain integer polynomials}

\author{
M. Idris\thanks{musbahu.idris@univ-lorraine.fr, IECL, Universit\'e de Lorraine, France}
\and
J.-M. Sac-\'Ep\'ee\thanks{jean-marc.sac-epee@univ-lorraine.fr, IECL, Universit\'e de Lorraine, France}
}

\date{}

\begin{document}

\maketitle

\begin{abstract}
Newman polynomials have coefficients in \(\{0,1\}\) and constant term \(1\),
whereas Littlewood polynomials have coefficients in \(\{-1,1\}\). We study
two questions concerning the existence of Newman and Littlewood multiples
for certain integer polynomials.

For Newman multiples, we revisit a question of Hare and
Mossinghoff~\cite{HareMossinghoff}: whether there exists a real number
\(\sigma>1\) such that every \(P\in\Z[x]\) with no nonnegative real root
and Mahler measure less than \(\sigma\) has a Newman multiple. To bound
any possible value of \(\sigma\) from above, we seek polynomials with no
nonnegative real root and no Newman multiple, and with Mahler measure as
small as possible. Using the updated database of known small-Mahler-measure
polynomials up to degree \(200\)~\cite{Known200Data}, we test, for each
listed polynomial \(p(x)\), both \(p(x)\) and its sign transform \(p(-x)\).
These two polynomials have the same Mahler measure, but the existence of a
Newman multiple is not invariant under the substitution \(x\mapsto -x\).
After a degree-bounded prefilter formulated as a mixed-integer linear
optimization problem, we apply the Hare--Mossinghoff certification
procedure to the remaining candidates. This yields \(14\) irreducible
polynomials of Mahler measure less than \(1.3\), with no nonnegative real
root and no Newman multiple. The smallest of their Mahler measures is
\(1.263095875491\ldots\), showing that any such \(\sigma\) is at most this
value.

For Littlewood multiples, we return to the three Newman polynomials listed
in Table~3 of Drungilas, Jankauskas, Junevi\v{c}ius, Klebonas and
\v{S}iurys~\cite{DrungilasEtAl2018}, for which the existence of a
Littlewood multiple of smallest possible degree was left unresolved. We
show that one of them attains this degree, whereas the other two have no
Littlewood multiple at either of their first two possible degrees.
\end{abstract}



\section{Introduction}

A Newman polynomial is a polynomial with coefficients in \(\{0,1\}\) and
constant term \(1\). A Littlewood polynomial is a polynomial all of whose
coefficients belong to \(\{-1,1\}\). We denote by \(\Ncal\) and \(\Lcal\)
the sets of Newman and Littlewood polynomials, respectively. If
\(P\in\Z[x]\), we say that \(P\) has a Newman multiple if it divides an
element of \(\Ncal\), and that it has a Littlewood multiple if it divides an
element of \(\Lcal\).

For
\[
  P(x)=a_d\prod_{j=1}^d(x-\alpha_j)\in\C[x],
\]
its Mahler measure is
\[
  M(P)=|a_d|\prod_{j=1}^d\max\{1,|\alpha_j|\}.
\]

The first part of the paper is motivated by a question of Hare and
Mossinghoff~\cite{HareMossinghoff}: is there a real number \(\sigma>1\) such
that every \(P\in\Z[x]\) with no nonnegative real root and
\(M(P)<\sigma\) has a Newman multiple? A natural way to obtain upper bounds
for any possible value of \(\sigma\) is to look for polynomials satisfying
the root condition and having no Newman multiple. If \(P\) is such a
polynomial, then the proposed property can hold only for
\(\sigma\leq M(P)\), so the relevant examples are those with \(M(P)\) as
small as possible.

Using Mossinghoff's updated database of known small-Mahler-measure
polynomials~\cite{Known200Data}, together with the unprinted sign transform
\(p(-x)\) of each listed polynomial \(p(x)\), we combine a degree-bounded
integer-linear prefilter with the Hare--Mossinghoff certification method. This yields \(14\) irreducible
polynomials with no nonnegative real root and no Newman multiple. The
smallest of their Mahler measures is \(1.263095875491\ldots\), showing that
any such \(\sigma\) is at most this value.

The second question addressed in this paper concerns the least possible
degree of a Littlewood multiple of a Newman polynomial. Drungilas,
Jankauskas and \v{S}iurys~\cite{DrungilasJankauskasSiurys} determined, for
every Newman polynomial of degree at most \(11\), whether it has a
Littlewood multiple. In later work, Drungilas, Jankauskas,
Junevi\v{c}ius, Klebonas and
\v{S}iurys~\cite{DrungilasEtAl2018} considered Newman polynomials of degree
at most \(10\) having Littlewood multiples and studied whether the smallest
possible degree of such a multiple is attained. They resolved all cases
except for the three polynomials listed in their Table~3. For these three
remaining cases, we show that one polynomial has a Littlewood multiple of
smallest possible degree, whereas the other two have no Littlewood multiple
at either of their first two possible degrees.

\section{A new upper bound in the Hare--Mossinghoff question}
\label{sec:hare-mossinghoff-bound}

\subsection{Background and motivation}\label{subsec:newman-background}

We first recall why small Mahler measure is relevant to the question of the existence of Newman multiples. Bombieri and
Vaaler~\cite{BombieriVaaler1987} proved that every \(f\in\Z[x]\) with
\(M(f)<2\) divides a nonzero polynomial with coefficients in
\(\{-1,0,1\}\). For Newman multiples, Hare and
Mossinghoff~\cite{HareMossinghoff} asked whether an analogous statement
holds below some threshold \(\sigma>1\), under the  assumption that
\(f\) has no nonnegative real roots.

The hypothesis that \(f\) have no nonnegative real root is necessary. Indeed, every Newman polynomial is strictly positive on \([0,\infty)\), so none of its divisors can vanish there. A second necessary condition comes from the root
location theorem of Odlyzko and Poonen. If
\[
 A_\rho=\{z\in\C:\rho^{-1}<|z|<\rho\},
 \qquad \tau=\frac{1+\sqrt5}{2},
\]
then every root of a Newman polynomial lies in the annulus
\(A_\tau\)~\cite{OdlyzkoPoonen}. Thus a polynomial with a Newman multiple must
have all of its roots in  \(A_\tau\setminus [0,\infty)\).

The appearance of the golden ratio \(\tau\) in this root-location result
suggests asking whether \(\tau\) can play, for Newman multiples, the same
role as the constant \(2\) in the height-one result of Bombieri and Vaaler.
Hare and Mossinghoff~\cite{HareMossinghoff} showed that \(\tau\) cannot
serve as such a threshold: they exhibited polynomials in
\(\Z[x]\), all of whose roots lie in the slit annulus
\(A_\tau\setminus [0,\infty)\), with Mahler measure less than \(\tau\), and
with no Newman multiple. The smallest Mahler measure in their table is
approximately \(1.556014485\). Drungilas, Jankauskas and
\v{S}iurys~\cite{DrungilasJankauskasSiurys} later found \(16\) Borwein
polynomials of degree at most \(9\) with no Newman multiple and with Mahler
measure below \(\tau\); the smallest measure among them is
\(1.436632261\ldots\). This was further lowered in~\cite{IdrisSacEpee} to
\(1.419404632\ldots\).

These examples make it natural to continue the search systematically among
polynomials of small Mahler measure. In~\cite{IdrisSacEpee}, we used the
small-Mahler-measure list then available~\cite{Known180Data}; it consisted
of \(8438\) known integer polynomials of degree at most \(180\)
and Mahler measure less than \(1.3\), all primitive, irreducible, and
noncyclotomic. Polynomials in the list having a nonnegative real root were
discarded at once, since such a root precludes the existence of a Newman
multiple. For every other polynomial in the list, with the exception of
three unresolved cases, we exhibited an explicit Newman multiple.

We now use the file \texttt{AllSmallMeasures\_200.txt} from Mossinghoff's
updated Mahler measure database~\cite{Known200Data}. It contains the known
primitive, irreducible, noncyclotomic integer polynomials of Mahler measure
less than \(1.3\) and degree at most \(200\). As in the previous version of
the list, only one representative of each pair \(\{p(x),p(-x)\}\) is
printed. The representative is chosen by the following sign convention: if
\(k\) is the least odd index for which the coefficient of \(x^k\) is
nonzero, then the listed polynomial is the one for which this coefficient
is positive~\cite{Known200Data}. Since \(p(x)\) and \(p(-x)\) have the same
Mahler measure but one may admit a Newman multiple while the other does
not, it is worthwhile to reconstruct and test the unprinted transforms as
well, thereby effectively doubling the pool of small-Mahler-measure
candidates. The \(14\) polynomials exhibited below all arise from this
sign-transformed family.

\subsection{A preliminary degree-bounded search}\label{subsec:newman-screening}
For the preliminary degree-bounded search, we restrict the augmented
family described above to polynomials of degree at most \(100\). We do not
consider polynomials of higher degree, since the Hare--Mossinghoff
certification procedure becomes difficult to apply in substantially larger
degrees. We begin by discarding any polynomial with a nonnegative real root. For each remaining polynomial \(f\), we search, up to a prescribed degree,
for a Newman polynomial divisible by \(f\), using the integer-linear
programming procedure described in~\cite[Section~2.1]{IdrisSacEpee}. This
preliminary search is used only as a filter: when the model is feasible, it
gives an explicit Newman multiple and the polynomial is discarded; when no
multiple is found in the tested range, no non-existence conclusion is drawn.

For a fixed target degree, the bounded search reduces to the following integer-linear feasibility test. Let \(f\) be one of the remaining
polynomials, and write \(n=\deg f\). For a prescribed degree \(D\geq n\),
put \(m=D-n\). We look for a multiplier
\[
  q(x)=x^m+b_{m-1}x^{m-1}+\cdots+b_1x+1\in\Z[x]
\]
such that
\[
  N(x)=f(x)q(x)
\]
is a Newman polynomial of degree \(D\).
Writing
\[
  f(x)=a_0+a_1x+\cdots+a_nx^n
  \qquad\text{and}\qquad
  N(x)=c_0+c_1x+\cdots+c_Dx^D,
\]
and comparing coefficients in the identity \(N=fq\), we see that each
coefficient \(c_j\) is obtained from the coefficients of \(f\) and the
coefficients of \(q\). Since the coefficients \(a_i\) are explicitly known,
each \(c_j\) is a linear expression in the integer unknowns
\(b_1,\ldots,b_{m-1}\). Thus, for this fixed value of \(D\), the search is
reduced to the feasibility of the linear system
\[
  0\leq c_j\leq 1\qquad(0\leq j\leq D),
  \qquad
  b_i\in\Z\qquad(1\leq i\leq m-1).
\]

In the preliminary filter used here, we ran this feasibility test for every
degree \(D\) with
\[
        \deg f\leq D\leq 250.
\]
If a Newman multiple was found in this range, the polynomial \(f\) was
discarded from the list of candidates. If no multiple was found, or if the
solver did not resolve one of the corresponding feasibility problems within
the time limit, the polynomial was kept for the certification step. This left
\(21\) candidates. No non-existence conclusion is drawn from this preliminary
filter; such a conclusion comes only from the Hare--Mossinghoff certification
method described next.

\subsection{The Hare--Mossinghoff certification method}\label{subsec:newman-stabilization}

We now turn to the second step. After the fast degree-bounded filter has
discarded all candidates for which it found an explicit Newman multiple,
we retain for certification only those remaining candidates having at
least one root outside the unit disk, as required by the
Hare--Mossinghoff method~\cite{HareMossinghoff}.

The Hare--Mossinghoff certification method is described in detail
in~\cite[Subsection~2.1]{HareMossinghoff}; we recall here only the
features needed for our computations. For each retained polynomial \(f\), 
let \(\beta\) denote a root with \(|\beta|>1\) chosen
for the computation, and put
\[
 I'(\beta)=
 \left\{z\in\C:|z|\leq\frac{|\beta|}{|\beta|-1}\right\}.
\]

For each nonnegative integer \(d\), let
\[
\Ncal'(\beta,d)
=
\{F(\beta):F\in\Ncal,\ \deg(F)\leq d\}\cap I'(\beta),
\]
and put
\[
\Ncal'(\beta)
=
\bigcup_{d\geq 0}\Ncal'(\beta,d).
\]
The sets \(\Ncal'(\beta,d)\) can be constructed recursively. Starting
with
\[
\Ncal'(\beta,0)=\{1\},
\]
we obtain
\[
\begin{aligned}
\Ncal'(\beta,d+1)
={}&
\Ncal'(\beta,d)\\
&{}\cup
\left(
\{\beta\omega:\omega\in\Ncal'(\beta,d)\}
\cap I'(\beta)
\right)\\
&{}\cup
\left(
\{\beta\omega+1:\omega\in\Ncal'(\beta,d)\}
\cap I'(\beta)
\right).
\end{aligned}
\]

The choice of \(I'(\beta)\) ensures that, if
\(\omega\notin I'(\beta)\), then neither \(\beta\omega\) nor
\(\beta\omega+1\) belongs to \(I'(\beta)\). 

As shown by Hare and Mossinghoff, the occurrence of \(0\) in the
iteration determines whether \(\beta\) is a root of a Newman polynomial.
Consequently, if, for some \(d\),
\[
 \Ncal'(\beta,d+1)=\Ncal'(\beta,d)
 \qquad\text{and}\qquad
 0\notin\Ncal'(\beta,d),
\]
then no Newman polynomial vanishes at \(\beta\). When \(f\) is irreducible
over \(\Z\), this certifies that \(f\) has no Newman multiple.

Among the \(21\) candidates left by the preliminary filter, we applied
the procedure, with maximum depth \(60\), to those having a root of
modulus greater than \(1\). It certified the \(14\) polynomials listed in
Table~\ref{tab:newman-polys}. No conclusion is drawn here for the other
candidates, so the table is not claimed to be exhaustive.

For every listed polynomial, the iteration stabilizes without reaching
\(0\). The column \(d_{\rm stab}\) gives the stabilization depth, and
\(|\Ncal'(\beta)|\) denotes the cardinality of the stabilized set.
Since all the polynomials listed in Table~\ref{tab:newman-polys} are
reciprocal, the order in which their coefficients are read is immaterial.
The entries are ordered by increasing Mahler measure.

\begingroup
\scriptsize
\setlength{\tabcolsep}{2pt}
\setlength{\LTleft}{0pt}
\setlength{\LTright}{0pt}
\begin{longtable}{@{}r r r r >{\ttfamily\raggedright\arraybackslash}p{0.47\linewidth}@{}}
\caption{Certified irreducible polynomials with no nonnegative real root and no Newman multiple.}\label{tab:newman-polys}\\
\toprule
\normalfont deg. & \normalfont Mahler measure &
\normalfont \(d_{\rm stab}\) & \normalfont \(|\Ncal'(\beta)|\) &
\normalfont coefficients \\
\midrule
\endfirsthead
\toprule
\normalfont deg. & \normalfont Mahler measure &
\normalfont \(d_{\rm stab}\) & \normalfont \(|\Ncal'(\beta)|\) &
\normalfont coefficients \\
\midrule
\endhead
\midrule
\multicolumn{5}{r}{\normalfont\footnotesize Continued on next page}\\
\endfoot
\bottomrule
\endlastfoot
44 & 1.263095875491 & 37 & 117962 & 1 0 -1 0 1 0 -1 0 0 0 0 -1 -1 1 1 -2 -1 2 1 -1 0 1 1 1 0 -1 1 2 -1 -2 1 1 -1 -1 0 0 0 0 -1 0 1 0 -1 0 1 \\
26 & 1.272019269348 & 20 & 1589 & 1 0 0 -1 0 -1 0 0 0 0 0 1 0 1 0 1 0 0 0 0 0 -1 0 -1 0 0 1 \\
50 & 1.273464959636 & 34 & 65657 & 1 -1 1 -1 0 0 0 0 0 0 -1 1 -1 1 -1 1 -1 1 0 0 1 -1 1 -1 1 -1 1 -1 1 -1 1 0 0 1 -1 1 -1 1 -1 1 -1 0 0 0 0 0 0 -1 1 -1 1 \\
48 & 1.279464310958 & 33 & 52354 & 1 -1 1 -1 0 0 0 0 0 0 -1 1 -1 1 -1 1 -1 1 0 0 1 -1 1 -1 1 -1 1 -1 1 0 0 1 -1 1 -1 1 -1 1 -1 0 0 0 0 0 0 -1 1 -1 1 \\
48 & 1.279702474008 & 31 & 34574 & 1 -1 1 -1 1 -1 0 0 -1 1 -2 2 -2 2 -1 1 0 0 1 -1 1 -1 1 -1 1 -1 1 -1 1 -1 1 0 0 1 -1 2 -2 2 -2 1 -1 0 0 -1 1 -1 1 -1 1 \\
44 & 1.284355030821 & 32 & 36221 & 1 0 -1 0 0 0 0 0 0 0 0 -1 0 0 0 1 0 0 0 0 0 0 1 0 0 0 0 0 0 1 0 0 0 -1 0 0 0 0 0 0 0 0 -1 0 1 \\
44 & 1.291273677148 & 31 & 28824 & 1 0 -1 0 0 0 0 0 0 -1 0 1 0 -1 0 1 0 0 1 0 -1 0 1 0 -1 0 1 0 0 1 0 -1 0 1 0 -1 0 0 0 0 0 0 -1 0 1 \\
56 & 1.294567098943 & 37 & 147782 & 1 0 -1 0 0 0 0 0 0 0 0 -1 0 0 0 1 0 0 0 0 0 0 1 0 0 0 0 0 -1 0 0 0 0 0 1 0 0 0 0 0 0 1 0 0 0 -1 0 0 0 0 0 0 0 0 -1 0 1 \\
48 & 1.295513327396 & 29 & 19823 & 1 0 -1 -1 0 1 1 0 -1 -1 0 1 0 0 0 0 0 0 0 0 0 0 0 0 1 0 0 0 0 0 0 0 0 0 0 0 0 1 0 -1 -1 0 1 1 0 -1 -1 0 1 \\
66 & 1.297441747412 & 39 & 210382 & 1 -1 0 0 0 0 0 0 0 0 0 0 0 0 0 0 -1 1 0 0 0 0 0 0 0 0 0 0 0 0 0 0 1 -1 1 0 0 0 0 0 0 0 0 0 0 0 0 0 0 1 -1 0 0 0 0 0 0 0 0 0 0 0 0 0 0 -1 1 \\
54 & 1.297657045330 & 30 & 22933 & 1 -1 0 0 0 0 0 0 0 -1 1 0 0 0 0 0 0 0 1 -1 0 0 0 0 0 0 1 -1 1 0 0 0 0 0 0 -1 1 0 0 0 0 0 0 0 1 -1 0 0 0 0 0 0 0 -1 1 \\
46 & 1.298706963798 & 28 & 14774 & 1 -1 1 -1 0 0 0 0 0 -1 0 0 0 1 0 0 0 0 0 1 0 0 0 -1 0 0 0 1 0 0 0 0 0 1 0 0 0 -1 0 0 0 0 0 -1 1 -1 1 \\
72 & 1.299208078270 & 40 & 263144 & 1 -1 0 0 0 0 0 0 0 0 0 0 -1 1 0 0 0 0 0 0 0 0 0 0 1 -1 0 0 0 0 0 0 0 0 0 1 -1 1 0 0 0 0 0 0 0 0 0 -1 1 0 0 0 0 0 0 0 0 0 0 1 -1 0 0 0 0 0 0 0 0 0 0 -1 1 \\
30 & 1.299764321006 & 18 & 811 & 1 -1 0 0 0 0 0 -1 1 0 0 0 0 0 1 -1 1 0 0 0 0 0 1 -1 0 0 0 0 0 -1 1 \\
\end{longtable}
\endgroup

The Hare--Mossinghoff certification procedure was implemented independently
in Julia and in GP/PARI 2.15.4. For each of the fourteen polynomials, both
implementations produced the same stabilization depth and the same
cardinality of the stabilized set.

The first entry of Table~\ref{tab:newman-polys} has Mahler measure
\(1.263095875491\ldots\). Since this polynomial has no nonnegative real root
and no Newman multiple, any value of \(\sigma\) for which the
Hare--Mossinghoff property holds must satisfy
\[
  \sigma\leq 1.263095875491\ldots .
\]

\section{Least-degree Littlewood multiples of Newman polynomials of degree
at most \(10\): the three cases that had remained unresolved}
\label{sec:littlewood}
\subsection{The modulo-\(2\) degree condition and previous results}
\label{subsec:littlewood-degree-condition}
Let \(\F_2\) denote the field with two elements. For \(P\in\Z[x]\), let
\(\widetilde P\in\F_2[x]\) denote its reduction modulo \(2\).

Let \(p\in\F_2[x]\) have nonzero constant term. We recall the
definition of \(\degtwo p\) introduced by Dubickas and
Jankauskas~\cite{DubickasJankauskas2009} and used by Drungilas,
Jankauskas, Junevi\v{c}ius, Klebonas and \v{S}iurys in their study of
the degree problem considered below
\cite[Section~4]{DrungilasEtAl2018}. Write the factorization of \(p\)
in \(\F_2[x]\) as
\[
  p(x)=(x+1)^m\prod_{j=1}^r\phi_j(x)^{m_j},
\]
where \(m\geq0\), the \(\phi_j\) are distinct irreducible polynomials
of degree at least \(2\), and \(m_j\geq1\). The product is empty if
\(r=0\). For each \(j\), let \(e_j\) be the unique odd positive
integer such that \(\phi_j\) divides the reduction of
\(\Phi_{e_j}\) modulo \(2\). Let \(s\) be the least nonnegative
integer such that
\[
  2^s\geq\max\{m+1,m_1,\ldots,m_r\},
\]
where the maximum is understood to be \(m+1\) when \(r=0\).
The quantity \(\degtwo p\) is then defined by
\[
  \degtwo p:=2^s\operatorname{lcm}(e_1,\ldots,e_r),
\]
with the usual convention that the least common multiple of an empty
family is \(1\).

The relevance of this quantity follows from
\cite[Lemma~12]{DubickasJankauskas2009}: for every integer \(n\geq0\),
\[
  p(x)\mid 1+x+\cdots+x^n
  \quad\Longleftrightarrow\quad
  \degtwo p\mid n+1
  \qquad\text{in }\F_2[x].
\]
If \(P\) divides a Littlewood polynomial \(L\), then
\(\widetilde P\) divides its reduction modulo \(2\), namely
\[
  1+x+\cdots+x^{\deg L}.
\]
Applying the preceding equivalence with \(p=\widetilde P\) gives
\[
  \degtwo\widetilde P\mid\deg L+1.
\]
Equivalently,
\begin{equation}\label{eq:degree-condition}
  \deg L=k\degtwo\widetilde P-1
  \qquad(k=1,2,\ldots).
\end{equation}
Thus \(\degtwo\widetilde P-1\) is the smallest possible degree of a
Littlewood multiple of \(P\), provided such a multiple exists.

Drungilas, Jankauskas and \v{S}iurys determined, for every Newman
polynomial of degree at most \(11\), whether it has a Littlewood
multiple~\cite{DrungilasJankauskasSiurys}. For a polynomial \(P\) of
degree \(d\), let
\[
  P^*(x)=x^dP(1/x)
\]
denote its reciprocal. A polynomial \(P\) has a Littlewood multiple of
degree \(D\) if and only if \(P^*\) does, so it is enough to consider one polynomial from each reciprocal pair.

Drungilas, Jankauskas, Junevi\v{c}ius, Klebonas and
\v{S}iurys~\cite{DrungilasEtAl2018} then investigated, for each Newman
polynomial \(P\) of degree at most \(10\) having a Littlewood multiple,
whether it has one of degree \(\degtwo\widetilde P-1\). With reciprocals omitted, they
settled all but the three cases listed in their
Table~3~\cite{DrungilasEtAl2018}. These three polynomials are reproduced
in Table~\ref{tab:three-cases}.

\begin{table}[ht]
\centering
\caption{The three cases from
\cite[Table~3]{DrungilasEtAl2018} for which the existence of a
Littlewood multiple of smallest possible degree was left open.}
\label{tab:three-cases}
\begin{tabular}{c p{0.52\linewidth} r r}
\toprule
Case & Polynomial \(P(x)\) & \(\degtwo\widetilde P\) &
\(\degtwo\widetilde P-1\) \\
\midrule
1 & \(x^9+x^8+x^6+x^5+x^4+x^3+x+1\) & 60 & 59 \\
2 & \(x^{10}+x^9+x^8+x^3+x^2+1\) & 1020 & 1019 \\
3 & \(x^{10}+x^9+x^7+x^6+x^5+1\) & 1022 & 1021 \\
\bottomrule
\end{tabular}
\end{table}

Let \(P_1,P_2,P_3\) denote the three polynomials in
Table~\ref{tab:three-cases}, in that order. We prove that \(P_1\) attains its
smallest possible degree, and that neither \(P_2\) nor \(P_3\) has a
Littlewood multiple at either of its first two possible degrees.

\subsection{A fixed-degree feasibility model}\label{subsec:littlewood-model}

Let \(P\) be a Newman polynomial, and let \(D\geq\deg P\) be fixed.
To determine whether \(P\) has a Littlewood multiple of degree \(D\), write
\[
  P(x)=\sum_{i\in S}x^i,
\]
where \(S\subseteq\{0,\ldots,\deg P\}\) is the set of indices of the
nonzero coefficients of \(P\), and seek
\[
  Q(x)=q_0+q_1x+\cdots+q_{D-\deg P}x^{D-\deg P}\in\Z[x]
\]
such that
\[
  L(x)=P(x)Q(x)=\sum_{j=0}^D \ell_jx^j
\]
is a Littlewood polynomial. We encode the coefficients
\(\ell_j\in\{-1,1\}\) by binary variables \(y_j\in\{0,1\}\), setting
\[
  \ell_j=2y_j-1.
\]
With the convention that \(q_t=0\) outside the range
\(0\leq t\leq D-\deg P\), comparison of coefficients gives
\begin{equation}\label{eq:milp-littlewood}
  \sum_{i\in S}q_{j-i}=2y_j-1
  \qquad(0\leq j\leq D).
\end{equation}

Thus the system consisting of equations
\eqref{eq:milp-littlewood} and the constraints
\[
  q_t\in\Z\quad(0\leq t\leq D-\deg P),
  \qquad
  y_j\in\{0,1\}\quad(0\leq j\leq D)
\]
is feasible if and only if \(P\) has a Littlewood multiple of degree
\(D\). By the symmetry \(L\mapsto-L\), we may impose the normalization
\(\ell_0=1\), or equivalently \(y_0=1\).

Whenever a feasible solution is found, we reconstruct \(L\) and \(Q\)
and verify the identity \(L=PQ\) by exact multiplication in \(\Z[x]\).

\subsection{The three cases}\label{subsec:littlewood-results}

\paragraph{The first polynomial.}
For \(P_1\), the system \eqref{eq:milp-littlewood} is feasible for
\(D=59\). This is certified by the following explicit identity. Let
\[
  L_1(x)=\sum_{j=0}^{59}\ell_jx^j,
  \qquad
  Q_1(x)=\sum_{j=0}^{50}q_jx^j.
\]
The coefficients are written in ascending order, with the sign chosen so
that \(\ell_0=1\):
\[
\begin{aligned}
(\ell_0,\ldots,\ell_{59})={}&(
1,1,1,1,1,1,1,1,1,1,\\
&1,-1,1,1,-1,1,-1,1,-1,-1,\\
&1,-1,1,-1,-1,1,-1,1,1,1,\\
&1,1,-1,1,1,1,1,-1,1,-1,\\
&1,-1,1,-1,-1,1,-1,1,1,1,\\
&1,1,1,1,1,1,1,1,-1,-1),
\end{aligned}
\]
and
\[
\begin{aligned}
(q_0,\ldots,q_{50})={}&(
1,0,1,-1,1,-2,2,-2,3,-3,\\
&4,-6,7,-7,8,-9,9,-10,9,-9,\\
&11,-11,12,-13,12,-13,14,-13,14,-13,\\
&13,-13,11,-10,9,-7,6,-5,3,-2,\\
&2,-1,1,-2,0,0,1,1,1,0,-1).
\end{aligned}
\]
Direct multiplication verifies that
\[
  L_1(x)=P_1(x)Q_1(x),
\]
and the displayed coefficients show that \(L_1\) is a Littlewood polynomial
of degree \(59\). Since \(59=\degtwo\widetilde P_1-1\), the least degree of
a Littlewood multiple of \(P_1\) is therefore \(59\).

\paragraph{The second and third polynomials.}
For \(P_2\) and \(P_3\), it is enough to apply the same fixed-degree
feasibility test at the second possible degree allowed by
\eqref{eq:degree-condition}, since non-existence at this degree also
rules out the first possible degree. Indeed, if \(N=\degtwo\widetilde P\) and \(L\)
is a Littlewood multiple of \(P\) of degree \(N-1\), then
\[
  (1+x^N)L(x)
\]
is again a Littlewood multiple of \(P\), now of degree \(2N-1\), since
the polynomials \(L(x)\) and \(x^NL(x)\) have disjoint supports. Thus
non-existence at degree \(2N-1\) also rules out degree \(N-1\).

We therefore solve the fixed-degree feasibility problem for \(P_2\) at
\[
  D=2\cdot1020-1=2039
\]
and for \(P_3\) at
\[
  D=2\cdot1022-1=2043.
\]

The corresponding integer feasibility problems were solved with
Gurobi. The branch-and-bound computations were carried to completion
and established that neither system is feasible.
Table~\ref{tab:gurobi-runs} gives the sizes of the two models and the
main computational data.

\begin{table}[ht]
\centering
\small
\caption{Computational data for the completed branch-and-bound
calculations for \(P_2\) and \(P_3\).}
\label{tab:gurobi-runs}
\begin{tabular}{c r r r r r}
\toprule
Polynomial & \(D\) & variables & nodes & simplex iterations & time (s) \\
\midrule
\(P_2\) & 2039 & 4070 & 62\,523 & 1\,426\,237 & 766.7 \\
\(P_3\) & 2043 & 4078 & 371\,892 & 8\,046\,184 & 1594.0 \\
\bottomrule
\end{tabular}
\end{table}

It follows that \(P_2\) has no Littlewood multiple of degree \(2039\), and
hence none of degree \(1019\); similarly, \(P_3\) has no Littlewood
multiple of degree \(2043\), and hence none of degree \(1021\).

\section{Computational details}\label{sec:implementation}

Most of the computations were carried out in Julia
1.11.6~\cite{Julia}. The Polynomials.jl package~\cite{PolynomialsJL}
was used for polynomial arithmetic and root computations, and the Nemo.jl
package~\cite{Nemo} for the irreducibility tests reported in
Table~\ref{tab:newman-polys}. The integer and mixed-integer linear models
were formulated with JuMP~\cite{JuMP} and solved with Gurobi
Optimizer~\cite{Gurobi}. A second, independent implementation of the
Hare--Mossinghoff procedure was written in GP/PARI
2.15.4~\cite{PARI}.

The time limit used in the preliminary Newman filter was \(180\) seconds
for each fixed-degree feasibility problem.

The calculations reported in Table~\ref{tab:gurobi-runs} were performed
with Gurobi Optimizer 12.0.3, with no time limit. For both runs, the
parameter settings were \texttt{Presolve=2}, \texttt{Cuts=2},
\texttt{MIPFocus=1}, and \texttt{NumericFocus=3}, with
\texttt{FeasibilityTol} and \texttt{IntFeasTol} both set to \(10^{-9}\).

All computations reported in this paper were performed on a dual-socket
server equipped with two Intel Xeon Gold 6138 processors at \(2.00\) GHz,
providing \(40\) physical cores in total, and \(204\) GiB of RAM.

\section{Conclusion}

The fourteen new counterexamples lower the current upper bound in the
question of Hare and Mossinghoff. A concrete continuation of this work is
to search for certified counterexamples of Mahler measure below
\(1.263095875491\ldots\). The main theoretical question remains whether
there exists a constant \(\sigma>1\) such that every integer polynomial
\(f\) with no nonnegative real root and \(M(f)<\sigma\) has a Newman
multiple.

For Littlewood multiples, the first degree permitted by reduction modulo
\(2\) need not be attained, and the computations for \(P_2\) and \(P_3\)
show that the second permitted degree need not be attained either. More
generally, if a Newman polynomial \(P\) has a Littlewood multiple, the
degree of every such multiple is of the form
\[
  k\degtwo\widetilde P-1,\qquad k\geq 1.
\]
The natural problem is therefore to determine the least value of \(k\)
for which such a multiple exists. For \(P_2\) and \(P_3\), the present
computations show that this value is greater than \(2\). Determining it
for these two polynomials, and understanding how large it can be in
general, remain open.

\end{document}